\documentclass[12pt]{article}
\usepackage{amsthm,amsmath,amssymb,amscd,amsfonts,mathrsfs,graphicx}
\usepackage[hidelinks]{hyperref}
\hypersetup{
     colorlinks   = true,
     citecolor    = blue
}
\hypersetup{linkcolor=blue}
\usepackage[all,cmtip]{xy}

\newtheorem{theorem}{Theorem}[section]
\newtheorem{definition}{Definition}
\newtheorem{lemma}[theorem]{Lemma}
\newtheorem{proposition}[theorem]{Proposition}
\newtheorem{corollary}[theorem]{Corollary}
\newcommand{\be}{\begin{equation}}
\newcommand{\ee}{\end{equation}}
\newcommand{\bes}{\begin{equation*}}
\newcommand{\ees }{\end{equation*}}
\newcommand{\bd}{\begin{definition}}
\newcommand{\ed}{\end{definition}}
\newcommand{\bp}{\begin{proposition}}
\newcommand{\ep}{\end{proposition}}
\newcommand{\bl}{\begin{lemma}}
\newcommand{\el}{\end{lemma}}
\newcommand{\bc}{\begin{corollary}}
\newcommand{\ec}{\end{corollary}}
\newcommand{\bt}{\begin{theorem}}
\newcommand{\et}{\end{theorem}}
\newcommand{\bpr}{\begin{proof}}
\newcommand{\epr}{\end{proof}}
\newcommand{\ben}{\begin{enumerate}}
\newcommand{\een}{\end{enumerate}}

\newcommand{\mb}{\mathbb}

\begin{document}
\setcounter{page}{1}
\title{\bf{Properties of Carath\'eodory measure hyperbolic universal covers of compact K\"ahler manifolds}}
\author{\it{Ngai-fung Ng}}
\date{}
\maketitle

\begin{abstract}
    This article explores some properties of universal covers of compact K\"ahler manifolds, under the assumption of Carath\'eodory measure hyperbolicity. In particular, by comparing invariant volume forms, an inequality is established between the volume of canonical bundle of a compact K\"ahler manifolds and the Carath\'eodory measure of its universal cover (similar result as in \cite{kikuta10}). Using similar method, an inequality is established between the restricted volume of canonical bundle of a compact K\"ahler manifolds and the restricted Carath\'eodory measure of its covering, solving a conjecture in \cite{kikuta13}.
\end{abstract}


\section{Introduction}
It is interesting to study non-compact complex manifolds with the assumption that there exists bounded holomorphic functions on it, since they include the examples of bounded domains and in particular, Hermitian symmetric spaces of non-compact type. In this article, we shall study properties of universal covers of compact K\"ahler manifolds, under the assumption of Carath\'eodory measure hyperbolicity. In this case, there are abundant supplies of holomorphic functions, leading to interesting properties. In \cite{kikuta10}, it was established that the volume of the canonical line bundle of a compact K\"ahler manifold is bounded from below by a constant multiple of its Carath\'eodory measure. It is natural to ask the same question regarding restricted volumes and measures. It was stated as a conjecture in \cite{kikuta13} that the restricted volume of the canonical line bundle of a compact K\"ahler manifold is also bounded from below by a constant multiple of its restricted Carath\'eodory measure. Inspired by the work of \cite{yeung}, which established uniform estimates among invariant metrics, the author obtains uniform estimates among invariant volumes forms and the conjecture follows as an easy consequence.

\section{Invariant Volume Forms}
Let $\mu^{d}_r$ be the \textbf{Poincare volume form} for $\mb{B}^d_r$, the complex $d$-ball of radius $r$, i.e.,
$$\mu^{d}_r= \frac{r^2}{(r^2-\|z\|^2)^{d+1}}  \cdot dz^1\wedge d\overline{z^1} \wedge\cdot\cdot\cdot \wedge dz^d\wedge d\overline{z^d}$$
For the unit ball $\mb{B}^d_1$, we write $\mu^{d}_1$ simply as $\mu^{d}$.\\\\
Let $M$ be an $n$-dimensional complex manifold. Here we introduce the following invariant volume forms:\\\\
The \textbf{Bergman pseudo-volume form} $v^B_{M}$ is defined by, for any $p\in M$,
$$ v^B_{M}(p):= \sum_i (\varphi_i \wedge \bar{\varphi_i})(p) $$
where $\{\varphi_i\}$ is an orthonormal basis for $L^2(M,K_M)$.\\\\
From \cite{hahn}, we have
$$v^B_M(p)= \varphi_p \wedge \overline{\varphi_p}$$ 
where $\varphi_p$ maximizes $(\varphi\wedge\bar{\varphi})(p)$ over the unit ball of $L^2(M, K_{M})$.\\\\
Suppose there exists a K\"ahler-Einstein metric $g^{KE}_M$ on $M$ satisfying 
$$ Ric(g^{KE}_{M}) = -2(n+1)\cdot \omega^{KE}_M $$ then the \textbf{K\"ahler-Einstein volumn form} $v^{KE}_{M}$ is defined by $$ v^{KE}_{M}:= \frac{(\omega^{KE}_M)^n}{n!} $$\\\\
The \textbf{Carath\'eodory pseudo-volume form} $v_{M}^C$ is defined by (\cite{eisenman} page 57), for any $p\in M$,
$$ v_{M}^C(p) := \sup \left\{ (f^* \mu^{n})(p); \hspace{2mm} f:M\to\mb{B}^n_1 \text{ holomorphic}, f(p)=0\right\} $$
$$\hspace{15mm} = \sup \left\{ |Jac(f)(p)|^2;  \hspace{2mm} f:M\to\mb{B}^n_1 \text{ holomorphic}, f(p)=0\right\} $$
where on any holomorphic coordinate system $\{z_1,...,z_n\}$ with $f=(f_1,...,f_n)$, we define $Jac(f)$ to be
$$
\begin{pmatrix}
\frac{\partial f_1}{\partial z_1} & \frac{\partial f_1}{\partial z_2} & ... & \frac{\partial f_1}{\partial z_n} \\
\frac{\partial f_2}{\partial z_1} & \frac{\partial f_2}{\partial z_2} & ... & \frac{\partial f_2}{\partial z_n} \\
... & ... & ... & ...\\
\frac{\partial f_n}{\partial z_1} & \frac{\partial f_n}{\partial z_2} & ... & \frac{\partial f_n}{\partial z_n}
\end{pmatrix}
$$ 
$|Jac(f)|$ to be 
$ \det(Jac(f))dz^1\wedge\cdot\cdot\cdot \wedge dz^n
$ and $|Jac(f)|^2$ to be
$$ |det(Jac(f))|^2  dz^1\wedge d\overline{z^1} \wedge\cdot\cdot\cdot \wedge dz^n\wedge d\overline{z^n}
$$
The \textbf{Kobayashi pseudo-volume form} $v_{M}^K$ is defined by (\cite{eisenman} page 57), for any $p\in M$,
$$ v_{M}^K(p) := \inf \left\{ |Jac(f)(0)|^{-2}; \hspace{2mm} f :\mb{B}^n_1 \to M \text{ holomorphic}, f(0)=p\right\} $$\\
\noindent Next let us consider the case when there is a $d$-dimensional subvariety $Z$ of $M$, in an analogous manner to the above, we may have the following definitions:\\\\
The \textbf{restricted Carath\'eodory pseudo-volume form} $v_{M|Z}^C$ on the regular locus $Z_{reg}$ of $Z$ is defined by, for any $p\in Z_{reg}$,
$$ v_{M|Z}^C(p) := \sup \left\{ |Jac(f|_{Z})(p)|^2;  \hspace{2mm} f:M\to\mb{B}^d_1 \text{ holomorphic}, f(p)=0\right\} $$
The \textbf{restricted Kobayashi pseudo-volume form} $v_{M|Z}^K$ on the regular locus $Z_{reg}$ of $Z$ is defined by, for any $p\in Z_{reg}$,
$$ v_{M|Z}^K(p) := \inf \left\{ |Jac(i_Z \circ f)(0)|^{-2}; \hspace{2mm}       f :\mb{B}^d_1 \to M \text{ holomorphic}, f(0)=p\right\} $$
where $i_Z:M\to Z$ is the restriction map.

\section{Basic Properties}
The following are the volume decreasing properties of the corresponding volume forms, which are included here for the convenience of the readers.
\bp[\bf Volume decreasing property for Carath\'eodory measure] Let $N$ be an $n$-dimensional complex manifold. Then, for all $\varphi\in Hol(M,N)$, $$ \varphi^* (v_{N}^C) \leq v_{M}^C$$
\ep
\bpr Any map from $N$ to $\mb{B}^n_1$ induces a map from $M$ to $\mb{B}^n_1$.
\epr

\bp[\bf Volume decreasing property for Kobayashi measure] Let $N$ be an $n$-dimensional complex manifold. Then, for all $\varphi\in Hol(M,N)$, $$  \varphi^* (v_{N}^K) \leq v_{M}^K$$
\ep
\bpr Any map from $\mb{B}^n_1$ to $M$ induces a map from $\mb{B}^n_1$ to $N$.
\epr

\noindent Similarly, we have the following:
\bp[\bf Volume decreasing property for restricted Carath\'eodory measure] Let $Y$ be an $n$-dimensional complex manifold and $W$ its $d$-dimensional complex subvariety. Then, for all $\varphi\in Hol(X,Y)$, such that $\varphi(Z)\subset W$, $$ (\varphi|_Z)^* (v_{Y|W}^C) \leq v_{X|Z}^C$$
\ep

\bp[\bf Volume decreasing property for restricted Kobayashi measure] Let $Y$ be an $n$-dimensional complex manifold and $W$ its $d$-dimensional complex subvariety. Then, for all $\varphi\in Hol(X,Y)$, such that $\varphi(Z)\subset W$, $$ (\varphi|_Z)^* (v_{Y|W}^K) \leq v_{X|Z}^K$$
\ep

\noindent We also have continuity property for the Carath\'eodory volume form. 

\bt $v^C_{M}$ is locally Lipschitz.
\et
\bpr Let $o\in M$, choose a local coordinate system $U$ at $o$. Let $\mb{B}^n_r$ be the coordinate ball centered at $o$. Let $p,q \in \mb{B}^n_{r/2}$ and suppose 
$$ v^C_{M}(p) \geq v^C_{M}(q) $$ then
\begin{equation*}
    \begin{split}
    v^C_{M}(p) - v^C_{M}(q) &= \sup{|Jac(f)(p)|^2} - \sup {\frac{|Jac(f)(q)|^2}{(1-\|f(q)\|^2)^{n+1}}} \\
    &\leq \sup\left\{ |Jac(f)(p)|^2 - |Jac(f)(q)|^2 \right\} 
    \end{split}
\end{equation*}where the supremum is taken over $\{f:M\to\mb{B}^n_1 \text{ holomorphic}, f(p)=0\}$.\\And for all $f:M\to\mb{B}^n_1$ holomorphic, we have
\begin{equation*}
    \begin{split}
    |Jac(f)(p)|^2 - |Jac(f)(q)|^2 
    & \leq \Big(|Jac(f)(p) - Jac(f)(q)| \Big) \Big(|Jac(f)(p)| + |Jac(f)(q)| \Big)\\
    &\leq \left( \frac{n^{1/2} 2^{n+2} n!}{r^{n+1}} \cdot \|p-q\| \right) \left( \frac{2^{n+1} n!}{r^{n}}  \right)\\
    &= \frac{n^{1/2}2^{2n+3}(n!)^2}{r^{2n+1}} \cdot \|p-q\| 
    \end{split}
\end{equation*}by Cauchy inequality applied to $|Jac(f)|$ and $f$.
\epr

\noindent Now let $X$ be an $n$-dimensional compact K\"ahler manifold, $\tilde{X}$ be its universal cover. The next result is a key property of Carath\'eodory measure hyperbolic universal cover.

\bp[\cite{kikuta10}] \label{caratheodory hyperbolic properties}Suppose $\tilde{X}$ is Carath\'eodory measure hyperbolic (i.e. $v^C_{\tilde{X}}>0$ at every point), then we have $c_1(K_X)>0$.
\ep
\bpr 
For any $f:\tilde{X}\to\mb{B}^n_1$ holomorphic, in a holomorphic coordinate system $\{z_1,...,z_n\}$, $f=(f_1,...,f_n)$ and let $\zeta\in\mb{C}^n$, we have
\begin{equation*}
\begin{split}
&\left(\sqrt{-1}\partial\bar\partial\log (f^*\mu^n)\right)(\zeta,\bar\zeta)  \\
=&(n+1)\sum_{i,j,k,l} \left(\frac{(1-\|f\|^2)\delta_{ij} + f_j\overline{f_i}}{(1-\|f\|^2)^2}\right)\frac{\partial f_i}{\partial z_k}\overline{\frac{\partial f_j}{\partial z_l}}\zeta^k \overline{\zeta^l} \\
=&(n+1)\left(\sum_{i,k,l} \frac{1}{1-\|f\|^2} \frac{\partial f_i}{\partial z_k}\overline{\frac{\partial f_i}{\partial z_l}}\zeta^k \overline{\zeta^l} + \frac{1}{(1-\|f\|^2)^2}\left| \sum_{i,k} \frac{\partial f_i}{\partial z_k} \overline{f_i} \zeta^k  \right|^2\right) \\
\geq& (n+1)\cdot\frac{1}{1-\|f\|^2}  \sum_{i,k,l} \frac{\partial f_i}{\partial z_k}\overline{\frac{\partial f_i}{\partial z_l}}\zeta^k \overline{\zeta^l} \\
=& (n+1)\cdot\frac{1}{1-\|f\|^2} \cdot \zeta^* (Jac(f))^*\cdot Jac(f) \zeta
\end{split}
\end{equation*}
By the Arzela-Ascoli theorem, for any $p\in\tilde{X}$, there is a map $f_p:\tilde{X}\to\mb{B}^n_1$ such that $v^C_{\tilde{X}}(p) = (f_p^*\mu^n)(p) = |Jac(f_p)|^2(p)$. From the above computation we see that locally $\log(f_p^*\mu^n)$ is strictly plurisubharmonic if and only if $|Jac(f_p)|\neq 0$.\\\\Let $p$ be in the local coordinate system and $C$ denote the Lipschitz constant as in the previous Theorem, using the estimates we obtained there, for $q \in B(p,\frac{v^C_{\tilde{X}}(p)}{10C}):=U$, 
$$|Jac(f_p)|^2(q) > |Jac(f_p)|^2(p) - \frac{v^C_{\tilde{X}}(p)}{10} = \frac{9}{10}\cdot v^C_{\tilde{X}}(p) $$
That is, on $U$, $|Jac(f_p)|^2 > \frac{9}{10}\cdot v^C_{\tilde{X}}(p)$.\\
For any $q, s$ in $B(p,\frac{v^C_{\tilde{X}}(p)}{10C})$, 
\begin{equation*}
\begin{split}
&\Big| |Jac(f_p)|^2(s) - |Jac(f_q)|^2(s) \Big| \\
\leq &\Big| |Jac(f_p)|^2(s) - |Jac(f_p)|^2(p) \Big| + \Big| |Jac(f_p)|^2(p) - |Jac(f_q)|^2(q) \Big| + \Big| |Jac(f_q)|^2(q) - |Jac(f_q)|^2(s) \Big| \\
\leq &\frac{v^C_{\tilde{X}}(p)}{10} + \frac{v^C_{\tilde{X}}(p)}{10} + \frac{v^C_{\tilde{X}}(p)}{5} \\
= &\frac{2}{5}\cdot v^C_{\tilde{X}}(p)
\end{split}
\end{equation*}
That is, for all $q\in U$, $|Jac(f_q)|^2>\frac{1}{2}\cdot v^C_{\tilde{X}}(p)>0$ on $U$.\\\\
Notice that $log (v^C_{\tilde{X}})|_U = \sup_{p\in U} log(f_p^*\mu^n)$ is continuous and therefore $log(v^C_{\tilde{X}})$ is strictly plurisubharmonic.\\
Interpret $v^C_{\tilde{X}}$ as a volume form over $X$ with $\log(v^C_{\tilde{X}})$ being strictly plurisubharmonic, that defines a strictly positive $(1,1)$ current on $K_X$ which can be approximated by a smooth volume form whose curvature form is positive (\cite{greeneandwu}, \cite{demailly92}). Hence, $c_1(K_X)>0$.
\epr

\noindent In the case that the covering is a bounded domain in $\mb{C}^n$, we have the following extra geometric property.

\bp[\textbf{Uniform Squeezing Property for $\tilde{X}$}] \label{uniformsqueezing} Suppose  $\tilde{X}$ is a bounded domain in $\mb{C}^n$ which covers a compact complex manifold $X$, then there exists $a,b$ satisfying $0 <a<b <\infty $ such that for any $x\in \tilde{X}$, there exists an embedding $\varphi_x: \tilde{X} \to \mathbb{C}^n$ with $\varphi_x(x)=0$ and $\mb{B}_a^n \subset \varphi_x(\tilde{X}) \subset \mb{B}_b^n$.
\ep
\bpr Let $A$ be any fundamental domain of $X$ in $\tilde{X}$. For any $x\in A$, take $r_x = \inf_{y \in \partial\tilde{X}}|x- y|>0$, $R_x = \sup_{y \in \partial\tilde{X}}|x-y|<\infty$ so that we have $\mb{B}^n_{r_x}(x)\subset \tilde{X}\subset \mb{B}^n_{R_x}(x)$.\\\\Now since $A$ is relatively compact in $\tilde{X}$, and $r_x, R_x$ are Lipschitz continuous in $x$, we have $\inf_A r_x:=a >0$ and $\sup_A R_x :=b <\infty$.\\\\For any $x \in \tilde{X}-A$, there is an automorphism of $\tilde{X}$ which brings $x$ to a point in $A$. And hence for all points in $\tilde{X}$, we have an embedding $\varphi_x$ such that $\varphi_x(x)=0$ and $\mb{B}^n_a \subset \varphi_x(\tilde{X}) \subset \mb{B}^n_b$.
\epr

\noindent Finally we include the statement of two Schwarz lemmas and a result on $L^2$-estimate for the $\bar\partial$-equation that will be used.
\bp[Schwarz lemma of \cite{mokyau}]
Let $M$ be a complete Hermitian manifold with scalar curvature bounded from below by $-K_1$ and let $N$ be a complex manifold of the same dimension with a volume form $V_N$ (i.e., positive $(n,n)$ form, $n=\dim N$) such that the Ricci form is negative definite and $(\frac{\sqrt{-1}}{2}\partial \bar\partial \log V_N)^n \geq K_2 V_N$. Suppose $f:M\to N$ is a holomorphic map and the Jacobian is nonvanishing at one point. Then $K_1>0$ and 
$$\sup \frac{f^*V_N}{V_M} \leq \frac{K_1^n}{n^nK_2} $$
\ep
\bp[Schwarz lemma of \cite{royden}]
Let $(M,g)$ be a complete K\"ahler manifold with Ricci curvature bounded from below by $k\leq 0$, and $(N,h)$ a K\"ahler manifold with holomorphic sectional curvature bounded from above by $K<0$. Then for any holomorphic map $f:M\to N$ we have
$$ \sum_{\alpha,\beta,i,j} g^{\alpha \bar{\beta}} f^i_{\alpha} \overline{f^j_{\beta}} h_{i\bar{j}} \leq \frac{2\nu}{\nu +1} \frac{k}{K} $$
where $\nu$ is the maximal rank of $df$.
\ep

\bp[\cite{hormander}, Theorem 4.1 of \cite{demailly82}]\label{l2} Let $M$ be a complete K\"ahler manifold and let $\varphi$ be a smooth strictly plurisubharmonic function on $M$. Then for any $v\in L^2_{n,1}(M,\partial\bar\partial\varphi,\varphi)$ with $\bar\partial v=0$, there is an $(n,0)$-form $u$ on $M$ such that $\bar\partial u =v$ and 
$$ \left|\int_M u\wedge \bar{u} e^{-\varphi} \right| \leq \int_M |v|^2_{\partial \bar\partial \varphi} e^{-\varphi}$$
\ep

\section{Computations on Complex N-ball}
Let $\mb{B}^n_r=\{z\in \mb{C}^n : \|z\| <r \}$. Define
$$\varphi = -\log(r^2-\|z\|^2)$$ and set 
$$\omega =\frac{\sqrt{-1}}{2} \partial\bar\partial \varphi = \frac{\sqrt{-1}}{2}\sum_{i,j}g_{i\bar{j}}  dz^i \wedge d\overline{z^j}$$
Now, Taylor series expansion gives
$$-\log(r^2-\|z\|^2)=-\log r^2 + \frac{\|z\|^2}{r^2} + \frac{\|z\|^4}{r^4}\cdot\frac{1}{2}+ \frac{\|z\|^6}{r^6}\cdot \frac{1}{3}+ \cdot\cdot\cdot + \frac{\|z\|^{2n}}{r^{2n}} \cdot\frac{1}{n}+\cdot\cdot\cdot $$
which shows that $$g_{i\bar{j}}(0) = \frac{\delta_{ij}}{r^2} $$ and $$dg(0)=0$$
The curvature tensor is
\bes
    \begin{split}
        R_{i\bar{j}k\bar{l}}(0) &= \frac{-\partial^4\varphi}{\partial z_i \partial \bar{z_j} \partial z_k \partial\bar{z_l}}(0) \\
        &=\frac{-1}{2r^4} \left(\frac{\partial^4 |z|^4}{\partial z_i \partial \bar{z_j} \partial z_k \partial\bar{z_l}}\right)(0) \\
        &=\frac{-1}{2r^4} \cdot
            \begin{cases}
            \hspace{2mm}4  \hspace{10mm}  i=j=k=l\\
            \hspace{2mm}2  \hspace{10mm}  \{i,k(\neq i)\}= \{j ,l\}\\
            \hspace{2mm}0  \hspace{10mm}  \text{otherwise}
            \end{cases}
    \end{split}
\ees
The holomorphic sectional curvature is 
$$\frac{4\cdot R(\partial_i,\partial_{\bar{i}},\partial_i,\partial_{\bar{i}})}{\|\partial_i+\partial_{\bar{i}}\|^4}(0)=\frac{-8/r^4}{4/r^4}=-2$$
The Ricci curvature tensor is
\bes
\begin{split}
Ric_{k\bar{l}}(0) &= g^{i\bar{j}}R_{i\bar{j}k\bar{l}}(0)\\
&= \sum_i r^2 \cdot R_{i\bar{i} k\bar{l}}(0)\\
&= 
\begin{cases}
0     \hspace{49mm} k\neq l \\
r^2\cdot \left[ \frac{-4}{2r^4} + \frac{-2}{2r^4} \cdot (n-1) \right] \hspace{9mm} k=l
\end{cases} \\
&= 
\begin{cases}
0     \hspace{49mm} k\neq l \\
\frac{-(n+1)}{r^2} \hspace{40mm} k=l
\end{cases}\\
&= 
\begin{cases}
0     \hspace{49mm} k\neq l \\
-(n+1) \cdot g_{k\bar{l}}(0)\hspace{22mm} k=l
\end{cases}
\end{split}
\ees
The scalar curvature is
$$ S(0)= \sum_{k,l} g^{k\bar{l}}(0) Ric_{k\bar{l}}(0)=\sum_k r^2\cdot  \frac{-(n+1)}{r^2}=-n(n+1)$$

\section{Uniform Estimates among Invariant Volumes}
In this section, we are going to establish estimates among various invariant volume forms.
\bt \label{uniform} Let $X$ be an $n$-dimensional compact K\"ahler manifold, $\tilde{X}$ be its universal cover. Then we have, at any $p\in \tilde{X}$, 
\ben
\item[(a)] \cite{kobayashi98} $$v_{\tilde{X}}^C \leq v_{\tilde{X}}^K$$
\item[(b)]
$$ v_{\tilde{X}}^C \leq \frac{1}{n!} \cdot v^{KE}_{\tilde{X}}$$
$$ v^{KE}_{\tilde{X}}\leq \frac{1}{n!}  \cdot v_{\tilde{X}}^K$$
$$ a_1 \cdot v_{\tilde{X}}^C \leq  v_{\tilde{X}}^B \leq a_2 \cdot v_{\tilde{X}}^C $$
if we assume $\tilde{X}$ is Carath\'eodory measure hyperbolic, where $a_1=a_1(X)>0$ and $a_2=a_2(X) <\infty$.
\item[(c)]
$$  v_{\tilde{X}}^K \leq \frac{b^{2n}}{a^{2n}} \cdot v_{\tilde{X}}^C$$
if we assume $\tilde{X}$ is a bounded domain in $\mb{C}^n$, where $a=a(X)>0,b=b(X)<\infty$ are as in Proposition (\ref{uniformsqueezing}).
\een
\et
\bpr
\noindent
\ben
\item[(a)]$v_{\tilde{X}}^C \leq v_{\tilde{X}}^K$:\\\\
Fix $\epsilon >0$, $p\in \tilde{X}$. Let $f:\mb{B}^n_1\to \tilde{X}$ and $g: \tilde{X} \to \mb{B}^n_1$ be any holomorphic maps with $f(0)=p$ and $g(p)=0$ respectively and such that they satisfy 
$$0 < |Jac(f)(0)|^{-2} < v_{\tilde{X}}^K + \epsilon$$ and $$|Jac(g)(p)|^2 > v_{\tilde{X}}^C - \epsilon$$
Consider the composition $ \mb{B}^n_1\overset{f}{\to} \tilde{X}\overset{g}{\to} \mb{B}^n_1$.\\
Applying Ahlfors-Schwarz lemma (Corollary 2.4.16 of \cite{kobayashi98}) gives
$$ (g\circ f)^*\mu^n \leq \mu^n$$
Expressed in terms of $f$ and $g$ gives
$$ |Jac(g)(p)|^2 \leq |Jac(f)(0)|^{-2}$$ and hence
$$ v_{\tilde{X}}^C - \epsilon \leq v_{\tilde{X}}^K + \epsilon $$
We are done since $\epsilon$ is arbitrary.\\

\item[(b)]$v_{\tilde{X}}^C \leq \frac{1}{n!} \cdot v^{KE}_{\tilde{X}}$:\\\\
By Proposition (\ref{caratheodory hyperbolic properties}), we know that $c_1(K_X)>0$.\\
By \cite{aubin}, \cite{yau}, there is a unique K\"ahler-Einstein metric $g^{KE}_X$ on $X$ such that $$ Ric(g^{KE}_X) = -2(n+1)\cdot\omega^{KE}_X$$
Pull back the metric by $\pi$ so that we have, on $\tilde{X}$, the complete metric $g^{KE}_{\tilde{X}}$ satisfying $$ Ric(g^{KE}_{\tilde{X}}) = -2(n+1)\cdot\omega^{KE}_{\tilde{X}}$$
Hence the scalar curvature of $\tilde{X}$ w.r.t this metric is $-n(n+1)$.\\\\
On $\mb{B}^n_1$, we can construct a K\"ahler-Einstein metric $g_{\mb{B}^n_1}^{KE}$ such that 
$$ Ric(g_{\mb{B}^n_1}^{KE}) = -2(n+1)\cdot \omega_{\mb{B}^n_1}^{KE}$$
Hence
$$ \left(\frac{\sqrt{-1}}{2}\partial\bar\partial log(\mu^n) \right)^n = (n+1)^n n! \cdot \mu^n $$
For any map $f:\tilde{X} \to \mb{B}^n_1$, applying the Schwarz lemma of \cite{mokyau} yields
$$ f^*\mu^n \leq \frac{(n(n+1))^n}{n^n (n+1)^n n!} \cdot v^{KE}_{\tilde{X}} = \frac{1}{n!}\cdot v^{KE}_{\tilde{X}}$$
Hence we have
$$ v_{\tilde{X}}^C \leq \frac{1}{n!} \cdot v^{KE}_{\tilde{X}}$$

\item[(b)]$v^{KE}_{\tilde{X}}\leq \frac{1}{n!}  \cdot v_{\tilde{X}}^K$:\\\\
For any $p\in \tilde{X}$, consider any map $f:\mb{B}^n_1\to \tilde{X}$ with $f(0)=p$. Applying the Schwarz lemma of \cite{mokyau} yields
$$ f^*v^{KE}_{\tilde{X}} \leq \frac{(n(n+1))^n}{n^n (n+1)^n n!}\cdot \mu^n$$
Hence we have
$$v^{KE}_{\tilde{X}}\leq \frac{1}{n!}  \cdot v_{\tilde{X}}^K$$
\item[(b)]$ a_1 \cdot v_{\tilde{X}}^C \leq  v_{\tilde{X}}^B \leq a_2 \cdot v_{\tilde{X}}^C $:\\\\
Recall that for any $p\in\tilde{X}$, there is an $f_p:\tilde{X}\to\mb{B}^n_1$ such that $f_p(p)=0$ and 
$v^C_{\tilde{X}}=|Jac(f_p)|^2>0 $. Let us first construct a bounded smooth strictly plurisubharmonic function $\psi$ on $\tilde{X}$. Pick any $p_1$ in $\tilde{X}$, $|Jac(f_{p_1})|\neq 0$ except perhaps on an ($n-1$)-dimensional subvariety $Z_1$. Pick a set of $p_{2i}$'s in the connected components of $Z_1$, then $|Jac(f_{p_1})|$ and $|Jac(f_{p_{2i}})|$'s all $\neq 0$ except perhaps on a strictly lower dimensional subvariety $Z_2$ than $Z_1$. Inductively, we will find a countable collection of $\{f_{q_1},...,f_{q_n},...\}$ such that their Jacobian determinant is not simultaneously zero. Consider the function $\psi := \log\Big(\sum_{i=1}^{\infty} \frac{1}{2^i}\|f_{q_i}\|^2 \Big)$. A straight forward computation shows that $\psi$ is plurisubharmonic and by our choice of the $f$'s, we know that $\psi$ is indeed strictly plurisubharmonic. Clearly $\psi$ is also bounded.\\\\
Let $\lambda$ be a $C^\infty$ real smooth cut off function that is $\equiv 1$ in a neighborhood of $p$. Choose $k>0$ such that $\varphi := k\psi + \lambda(n+1)\log|z|^2$ is strictly psh on $\tilde{X}\backslash\{p\}$. By Proposition (\ref{l2}), we can solve the equation $\bar\partial u = \bar\partial (\lambda |Jac(f_p)|)$ on $\tilde{X}\backslash\{p\}$ with estimate
$$\left| \int_{\tilde{X}\backslash\{p\}}u\wedge \bar{u} e^{-\varphi}\right| \leq \int_{\tilde{X}\backslash\{p\}} |\bar\partial (\lambda |Jac(f_p)|)|^2_{\partial \bar\partial \varphi} e^{-\varphi} \leq C$$ 
This implies that the $L^2$ holomorphic $n$-form $v = \lambda |Jac(f_p)| - u$ satisfies $$v(p) = |Jac(f_p)|(p)\neq 0$$ because from the choice of $\varphi$ we know that $u(p)=0$. Hence, $v^B_{\tilde{X}}(p)>0$. Consider the ratio $v^B_{\tilde{X}}/v^C_{\tilde{X}}$, being continuous and invariant under automorphisms of $\tilde{X}$, let $A$ be a fundamental domain of $X$ in $\tilde{X}$, we have
$$ \inf_{\tilde{X}} v^B_{\tilde{X}}/v^C_{\tilde{X}} = \inf_{A} v^B_{\tilde{X}}/v^C_{\tilde{X}} := a_1(X) > 0$$ and 
$$ \sup_{\tilde{X}} v^B_{\tilde{X}}/v^C_{\tilde{X}} = \sup_{A} v^B_{\tilde{X}}/v^C_{\tilde{X}} := a_2(X) < \infty$$

\item[(c)]$v_{\tilde{X}}^K \leq \frac{b^{2n}}{a^{2n}} \cdot v_{\tilde{X}}^C$:\\\\
Realize $\tilde{X}$ as $\mb{B}_a^n \subset \tilde{X}\subset \mb{B}_b^n$ with $p=(0,...,0)$. Consider the composition
$$ \mb{B}^n_1\to \mb{B}^n_a \hookrightarrow \tilde{X}\hookrightarrow \mb{B}^n_b \to \mb{B}^n_1$$ in which the first and last maps are given, respectively, by
$$ (w_1,...,w_n) \mapsto (aw_1,...,aw_n) $$ and 
$$ (z_1,...,z_n) \mapsto (z_1/b,...,z_n/b) $$
Hence, 
$$ v_{\tilde{X}}^C(p) \geq \frac{1}{b^{2n}} \cdot dz^1\wedge d\overline{z^1}\wedge\cdot\cdot\cdot\wedge dz^n\wedge d\overline{z^n}$$ and
$$ v_{\tilde{X}}^K(p) \leq \frac{1}{a^{2n}} \cdot dz^1\wedge d\overline{z^1}\wedge\cdot\cdot\cdot\wedge dz^n\wedge d\overline{z^n}$$

\een
\epr

\noindent We know that the Carath\'eodory pseudo-volume form is invariant under biholomorphic map. Hence, the \textbf{Carath\'eodory measure} $\mu_{\tilde{X}}^C(X)$ of $X$ can thus be defined to be
$$ \mu_{\tilde{X}}^C(X):= \int_A v_{\tilde{X}}^C $$
where $A$ can be any fundamental domain of $X$ in $\tilde{X}$.\\\\
\noindent Next, let $L$ be a holomorphic line bundle over $X$, the \textbf{volume} $vol_{X}(L)$ of $L$ is defined as 
$$ vol_{X}(L):=\limsup_{m\to \infty}\frac{\dim H^0(X, \mathcal{O}(mL))}{m^n/n!}$$
If $L>0$, by the Hirzebruch-Riemann-Roch formula and Kodaira vanishing theorem, we have 
$$ vol_{X}(L) = \int_X c_1(L)^n $$


\bc
Let $X$ be an $n$-dimensional compact K\"ahler manifold such that $\tilde{X}$ is Carath\'eodory measure hyperbolic, then we have
$$ \frac{(n!)^2(n+1)^n}{(\pi)^n}\cdot \mu_{\tilde{X}}^C(X) \leq vol_{X}(K_X)$$
\ec
\bpr 
From 
$$n! \cdot v_{\tilde{X}}^C \leq v^{KE}_{\tilde{X}}= \frac{1}{n!} (\omega^{KE}_{\tilde{X}})^n =  \frac{1}{n!(2(n+1))^n} \left(-Ric(g^{KE}_{\tilde{X}}) \right)^n$$ we have
$$ \frac{(n!)^2(n+1)^n}{(\pi)^n}\cdot v_{\tilde{X}}^C \leq \left(-\frac{1}{2\pi} Ric(g^{KE}_{\tilde{X}}) \right)^n $$
Let $A$ be a fundamental domain of $X$ in $\tilde{X}$. Integrate the above inequality over $A$, we have 
$$ \frac{(n!)^2(n+1)^n}{(\pi)^n} \int_A v_{\tilde{X}}^C \leq \int_{A}\left(-\frac{1}{2\pi} Ric(g^{KE}_{\tilde{X}}) \right)^n$$
Since $c_1(K_X)>0$, we get
$$ \frac{(n!)^2(n+1)^n}{(\pi)^n}\cdot \mu_{\tilde{X}}^C(X) \leq vol_{X}(K_X) $$
\epr
\noindent Similar estimate is obtained in \cite{kikuta10} by using a different method.

\section{Uniform Estimates among Restricted Invariant Volumes}
Similar to the treatment of the previous section, we have the following:
\bt \label{uniformrestricted} Let $X$ be an $n$-dimensional compact K\"ahler manifold, $\pi:\tilde{X}\to X$ be its universal cover, $Z$ be a $d$-dimensional subvariety of $X$, $\tilde{Z}_{reg}$ is the regular part of the subvariety $\tilde{Z}:=\pi^{-1}(Z)$. Then we have, at any $p\in \tilde{Z}_{reg}$, 

\ben
\item[(a)] $$ v_{\tilde{X}|\tilde{Z}}^C \leq v_{\tilde{X}|\tilde{Z}}^K $$
\item[(b)]
$$v_{\tilde{X}|\tilde{Z}}^C \leq \frac{d^d(n+1)^d}{(d+1)^d}\cdot v^{KE}_{{\tilde{X}}|{\tilde{Z}}}$$
if we assume $\tilde{X}$ is Carath\'eodory measure hyperbolic.
\item[(c)]
$$ v_{\tilde{X}|\tilde{Z}}^K \leq \frac{b^{2d}}{a^{2d}}\cdot v_{\tilde{X}|\tilde{Z}}^C$$
if we assume $\tilde{X}$ is a bounded domain in $\mb{C}^n$, where $a=a(X),b=b(X)$ are as in Proposition (\ref{uniformsqueezing}).
\een
\et
\bpr
\noindent
\ben
\item[(a)]$v_{\tilde{X}|\tilde{Z}}^C \leq v_{\tilde{X}|\tilde{Z}}^K$:\\\\
Fix $\epsilon >0$, $p\in \tilde{Z}_{reg}$. Let $f:\mb{B}^d_1\to \tilde{X}$ and $g: \tilde{X} \to \mb{B}^d_1$ be any holomorphic maps with $f(0)=p$ and $g(p)=0$ respectively and such that they satisfy 
$$0 < |Jac(i_{\tilde{Z}_{reg}} \circ f)(0)|^{-2} < v_{\tilde{X}|\tilde{Z}}^K + \epsilon$$ and $$|Jac(g|_{\tilde{Z}_{reg}})(p)|^2 > v_{\tilde{X}|\tilde{Z}}^C - \epsilon$$
Consider the composition $ \mb{B}^d_1 \xrightarrow{f} \tilde{X}\xrightarrow{i_{\tilde{Z}_{reg}}} \tilde{Z}_{reg} \xrightarrow{g|_{\tilde{Z}_{reg}}} \mb{B}^d_1$.\\
Applying Ahlfors-Schwarz lemma (Corollary 2.4.16 of \cite{kobayashi98}) gives
$$ (g\circ f)^*\mu^d \leq \mu^d$$
Expressed in terms of $i_{\tilde{Z}_{reg}} \circ f$ and $g|_{\tilde{Z}_{reg}}$ gives
$$ |Jac(g|_{\tilde{Z}_{reg}})(p)|^2 \leq |Jac(i_{\tilde{Z}_{reg}} \circ f)(0)|^{-2}$$ and hence
$$ v_{\tilde{X}|\tilde{Z}}^C - \epsilon \leq v_{\tilde{X}|\tilde{Z}}^K + \epsilon $$
We are done since $\epsilon$ is arbitrary.\\

\item[(b)]$v_{\tilde{X}|\tilde{Z}}^C \leq \frac{d^d(n+1)^d}{(d+1)^d}\cdot v^{KE}_{{\tilde{X}}|{\tilde{Z}}}$:\\\\
By Proposition (\ref{caratheodory hyperbolic properties}), we know that $c_1(K_X)>0$.\\
By \cite{aubin}, \cite{yau}, there is a unique K\"ahler-Einstein metric $g^{KE}_X$ on $X$ such that $$ Ric(g^{KE}_X) = -2(n+1)\cdot\omega^{KE}_X$$
Pull back the metric by $\pi$ so that we have, on $\tilde{X}$, the complete metric $g^{KE}_{\tilde{X}}$ satisfying $$ Ric(g^{KE}_{\tilde{X}}) = -2(n+1)\cdot\omega^{KE}_{\tilde{X}}$$\\
On $\mb{B}^d_1$, we can construct a K\"ahler-Einstein metric $g_{\mb{B}^d_1}^{KE}$ such that 
$$ Ric(g_{\mb{B}^d_1}^{KE}) = -2(d+1)\cdot \omega_{\mb{B}^d_1}^{KE}$$ \\
For any map $f:\tilde{X} \to \mb{B}^d_1$, applying the Schwarz lemma of \cite{royden} yields
$$ f^*\omega_{\mb{B}^d_1}^{KE} \leq \frac{2d}{d+1} \frac{n+1}{2} \cdot \omega^{KE}_{\tilde{X}} $$
Restricting to $\tilde{Z}_{reg}$ and taking $d$-th power on both sides yields
$$ f^* \mu^d \leq \frac{d^d(n+1)^d}{(d+1)^d}\cdot \frac{(\omega^{KE}_{\tilde{X}})^d}{d!}$$ 
and therefore
$$v_{\tilde{X}|\tilde{Z}}^C \leq \frac{d^d(n+1)^d}{(d+1)^d}\cdot v^{KE}_{{\tilde{X}}|{\tilde{Z}}}$$

\item[(c)]$ v_{\tilde{X}|\tilde{Z}}^K \leq \frac{b^{2d}}{a^{2d}}\cdot v_{\tilde{X}|\tilde{Z}}^C$:\\\\
Realize $\tilde{X}$ as $\mb{B}_a^n \subset \tilde{X}\subset \mb{B}_b^n$ with $p=(0,...,0)$ and locally at $p$, $\tilde{Z}$ is given by $z_{d+1}=...=z_n=0$. Consider the composition
$$ \mb{B}^d_1\to \mb{B}^n_a \hookrightarrow \tilde{X}\hookrightarrow \mb{B}^n_b \to \mb{B}^d_1$$ in which the first and last maps are given, respectively, by
$$ (w_1,...,w_d) \mapsto (aw_1,...,aw_d,0,...,0) $$ and 
$$ (z_1,...,z_n) \mapsto (z_1/b,...,z_d/b) $$
Hence, 
$$v_{\tilde{X}|\tilde{Z}}^C(p) \geq \frac{1}{b^{2d}} \cdot dz^1\wedge d\overline{z^1}\wedge\cdot\cdot\cdot\wedge dz^n\wedge d\overline{z^n}$$ and
$$v_{\tilde{X}|\tilde{Z}}^K(p) \leq \frac{1}{a^{2d}} \cdot dz^1\wedge d\overline{z^1}\wedge\cdot\cdot\cdot\wedge dz^n\wedge d\overline{z^n}$$
\een
\epr

\noindent We know that the restricted Carath\'eodory pseudo-volume form is invariant under biholomorphic map preserving $Z$. Hence, the \textbf{restricted Carath\'eodory measure} $\mu_{\tilde{X}|\tilde{Z}}^C(Z)$ of $Z$ can thus be defined to be 
$$ \mu_{\tilde{X}|\tilde{Z}}^C(Z):= \int_{A\cap \tilde{Z}_{reg}} v_{\tilde{X}|\tilde{Z}}^C$$
where $A$ can be any fundamental domain of $X$ in the covering $\pi:\tilde{X}\to X$ and $\tilde{Z}_{reg}$ is the regular part of the subvariety $\tilde{Z}:=\pi^{-1}(Z)$.\\\\
Next, let $L$ be a holomorphic line bundle over $X$, the \textbf{restricted volume} $vol_{X|Z}(L)$ of $L$ along $Z$ is defined as 
$$ vol_{X|Z}(L):=\limsup_{m\to \infty}\frac{\dim H^0(X|Z, \mathcal{O}(mL))}{m^d/d!}$$
$$:=\limsup_{m\to \infty}\frac{\dim Im[i^*_Z:H^0(X, mL) \to H^0(Z, mL|_Z)]}{m^d/d!}$$
The following result will be needed to estimate the restricted volume of $K_X$.

\bp[\cite{boucksom},\cite{hisamoto},\cite{matsumura}] \label{restrictedvolume}Let $X$ be an $n$-dimensional projective manifold, $L$ a big line bundle over $X$ (i.e. its Kodaira dimension is $n$), and $Z$ an irreducible closed $d$-dimensional subvariety of $X$. Furthermore, suppose that $Z\not\subset \mb{B}_+(L)$. Then
$$ vol_{X|Z}(L)=\sup_T\left\{ \int_{Z_{reg}} (T|_{Z_{reg}})_{ac}^d \right\}$$
where the supremum is taken over all $T$'s that are semi-positive $(1,1)$-current in $c_1(L)$ that have measure zero unbounded loci and whose unbounded loci do not contain $Z$. Here we denote $T|_{Z_{reg}}$ the restriction of $T$ to the regular locus of $Z$ and $(T|_{Z_{reg}})_{ac}$ its absolutely continuous part.
\ep

\noindent Here we arrive at our promised result, settling a conjecture in \cite{kikuta13}.

\bc Let $X$ be an $n$-dimensional compact K\"ahler manifold such that $\tilde{X}$ is Carath\'eodory measure hyperbolic, $Z$ be a $d$-dimensional subvariety of $X$, then we have
$$  \frac{d!(d+1)^d}{\pi^d d^d} \cdot \mu^C_{\tilde{X}|\tilde{Z}}(Z) \leq vol_{X|Z}(K_X) $$
\ec

\bpr
From
$$ \frac{(d+1)^d}{d^d(n+1)^d} \cdot v_{\tilde{X}|\tilde{Z}}^C \leq v^{KE}_{{\tilde{X}}|{\tilde{Z}}} = \frac{1}{d!} ((\omega^{KE}_{\tilde{X}})|_{\tilde{Z}_{reg}})^d = \frac{1}{d!(2(n+1))^d}\left(-Ric(g^{KE}_{\tilde{X}})|_{\tilde{Z}_{reg}}\right)^d$$
we have
$$ \frac{d!(d+1)^d}{\pi^d d^d} \cdot v_{\tilde{X}|\tilde{Z}}^C \leq \left( -\frac{1}{2\pi}Ric(g^{KE}_{\tilde{X}})|_{\tilde{Z}_{reg}} \right)^d$$ 
Let $A$ be a fundamental domain of $X$ in $\tilde{X}$.
Integrate the above inequality over $A\cap \tilde{Z}_{reg}$ gives 
$$ \frac{d!(d+1)^d}{\pi^d d^d} \int_{A\cap \tilde{Z}_{reg}} v_{\tilde{X}|\tilde{Z}}^C \leq \int_{A\cap \tilde{Z}_{reg}} \left( -\frac{1}{2\pi}Ric(g^{KE}_{\tilde{X}})|_{\tilde{Z}} \right)^d $$
Apply Proposition (\ref{restrictedvolume}), we get
$$  \frac{d!(d+1)^d}{\pi^d d^d} \cdot \mu^C_{\tilde{X}|\tilde{Z}}(Z) \leq vol_{X|Z}(K_X) $$
\epr

\section{Completeness of Bergman Metric}
Suppose that $v^B_{M}>0$ at every point, the \textbf{Bergman pseudometric} $g^B_{M}$ on $M$ is defined by
$$ (g^B_{M})_{i\bar{j}} := \frac{\partial^2 \log
(v^B_{M})}{\partial z_{i} \partial \overline{z_j} } $$
From \cite{hahn}, we have
$$g^B_{M}(p,v)= \frac{\partial_v g_{p,v} \wedge \overline{\partial_v g_{p,v}}}{\varphi_p \wedge \overline{\varphi_p}}$$ 
where $\varphi_p$ maximizes $(\varphi\wedge\bar\varphi)(p)$ and $g_{p,v}$ maximizes $(\partial_v g \wedge \overline{\partial_v g})(p)$ with $g(p)=0$ over the unit ball of $L^2(M, K_{M})$.

\bt
Let $X$ be an $n$-dimensional compact K\"ahler manifold such that $\tilde{X}$ is Carath\'eodory measure hyperbolic, then $g^B_{\tilde{X}}$ is complete.
\et 
\bpr By Theorem (\ref{uniform}), we know that $v^B_{\tilde{X}}>0$ at every point, so $g^B_{\tilde{X}}$ is well-defined. By Proposition (\ref{caratheodory hyperbolic properties}), we know that $c_1(K_X)>0$, hence there exists a complete K\"ahler-Einstein metric on $\tilde{X}$. By Proposition 2.3 of \cite{chen}, to show that the Bergman metric is positive definite, it suffices to construct a bounded smooth strictly plurisubharmonic function $\psi$ on $\tilde{X}$. As is discussed in the proof of Theorem (\ref{uniform}), such a function $\psi$ exists.\\\\
Now let $\{y_i\}$ be a Cauchy sequence in $\tilde{X}$ w.r.t. $g^B_{\tilde{X}}$. As $g^B_{\tilde{X}}$ is an invariant metric, the push-forward of $g^B_{\tilde{X}}$ by $\pi:\tilde{X}\to X$ onto $X$ is a well-defined, positive definite metric, denoted by $g_X$. Hence, $\{\pi(y_i)\}$ is Cauchy w.r.t. $g_X$ on the compact manifold $X$. Hence, $\{\pi(y_i)\}$ converges to, say, $z\in X$. Let $\pi^{-1}\{z\}$ = $\{x_i\}_{i\in I}$. By the discreteness of the action of deck transformations on $\{x_i\}_{i\in I}$, there exists $\epsilon$-neigborhoods of the $x_i$'s that are mutually non-intersecting, hence, there exists $x_{i_0}$ in $\{x_i\}_{i\in I}$ that is the limit of $\{y_i\}$.
\epr 

\section{Invariant metrics}
For a tangent vector $v\in T_pM$ on a complex manifold $M$, the \textbf{Kobayashi and Carath\'eodory pseudometrics} are defined by
$$ \sqrt{g^C_M(p,v)} =\sup \{\|u\| ; \hspace{2mm} f:M\to\mb{B}^1_1 \text{ holomorphic},  f(p)=0, f_*v=u\} $$
$$ \sqrt{g^K_M(p,v)} =\inf \{ \|u\|; \hspace{2mm}  f:\mb{B}^1_1 \to M \text{ holomorphic}, f(0)=p, f_*u=v\} $$\\
\noindent The following theorem includes various estimates established among invariant metrics, they are included for the convenience of the readers. The only part that is new is part(b), which is proved by a simple application of Schwarz lemma of \cite{royden}.
\bt \label{quasiisometric} Let $X$ be an $n$-dimensional compact K\"ahler manifold, $\tilde{X}$ be its universal cover. Then we have, at any $p\in \tilde{X}$ and $v \in T_p\tilde{X}$, 
\ben
\item[(a)] $$ g^C_{\tilde{X}} \leq g^K_{\tilde{X}} $$
$$ g^C_{\tilde{X}} \leq g^B_{\tilde{X}} $$
\item[(b)]
$$ g^C_{\tilde{X}} \leq (n+1) \cdot g^{KE}_{\tilde{X}} $$
if we assume $\tilde{X}$ is Carath\'eodory measure hyperbolic.
\item[(c)]
$$  g^K_{\tilde{X}} \leq \frac{b^2}{a^2} \cdot g^C_{\tilde{X}} $$
$$ g^B_{\tilde{X}} \leq \left[\frac{2\pi}{a^3}\left(\frac{2b}{a}\right)^n\right]^2 g^K_{\tilde{X}} $$
$$ \frac{a^2}{b^2n}g^K_{\tilde{X}}\leq g^{KE}_{\tilde{X}}  \leq \frac{b^{4n-2}n^{n-1}}{a^{2n-2}}g^K_{\tilde{X}}$$
if we assume $\tilde{X}$ is a bounded domain in $\mb{C}^n$, where $a=a(X),b=b(X)$ are as in Proposition (\ref{uniformsqueezing}).
\een
\et
\bpr
\noindent
\ben
\item[(a)]$ g^C_{\tilde{X}} \leq g^K_{\tilde{X}} $:\\\\It follows directly from the Schwarz lemma.
\item[(a)]$ g^C_{\tilde{X}} \leq g^B_{\tilde{X}} $:\\\\See \cite{hahn}.
\item[(b)]$ g^C_{\tilde{X}} \leq (n+1) \cdot g^{KE}_{\tilde{X}} $:\\\\For any map $f:\tilde{X} \to \mb{B}^1_1$, applying the Schwarz lemma of \cite{royden} yields
$$ f^*\mu^1 \leq \frac{2(1)(-(2n+2))}{(1+1)(-2)} \cdot g^{KE}_{\tilde{X}}$$
Hence we have
$$ g^C_{\tilde{X}} \leq (n+1) \cdot g^{KE}_{\tilde{X}}$$
\item[(c)]
See \cite{yeung}.
\een
\epr

\end{document}